\documentclass[11pt]{article}
\usepackage{amsmath,amsfonts,amsthm,graphicx,amssymb}
\newtheorem{thm}{Theorem}[section]

\newtheorem{lem}[thm]{Lemma}
\newtheorem{prop}[thm]{Proposition}
\newtheorem{defn}[thm]{Definition}

\newtheorem{exap}[thm]{Example}
\newtheorem{rem}[thm]{Remark}
\numberwithin{equation}{section}

\newcommand\bP{{\mathbb P}}
\newcommand\Eff{{\rm Pseff}}
\newcommand\Nef{{\rm {Nef}}}
\newcommand\Neff{{\rm {Nef}}}
\newcommand\rk{{\rm {rank}}}
\newcommand\Proj{{\rm {Proj}}}
\newcommand\Sym{{\rm {Sym}}}
\newcommand\pr{{\rm {pr}}}
\newcommand\cone{{\rm {cone}}}

\input xy
\xyoption{all}
\date{}
\title{The cones of effective cycles on projective bundles over curves}
\author{Mihai Fulger}
\begin{document}
\maketitle

\section{Introduction}
\par The study of the cones of curves or divisors on complete varieties is a classical subject in Algebraic Geometry (cf. \cite{Kleiman}, \cite{Kawamata}, \cite{CP}) and it still is an active research topic  (cf. \cite{BDPP}, \cite{Lazarsfeld1} or \cite{Bou}). However, little is known if we pass to higher (co)dimension. In this paper we study this problem in the case of projective bundles over curves and describe the cones of effective cycles in terms of the numerical data appearing in a Harder--Narasimhan filtration. This generalizes to higher codimension results of Miyaoka and others (\cite{Miyaoka}, \cite{Bruzzo}) for the case of divisors. An application to projective bundles over a smooth base of arbitrary dimension is also given.

\par Given a smooth complex projective variety $X$ of dimension $n$, consider the real vector spaces 
$$N^k(X):=
\frac{\mbox{The real span of }\{[Y]\ :\ Y\mbox{ is a subvariety of } X\mbox{ of codimension }k\}}
{\mbox{numerical equivalence of cycles}}.$$ We also denote it by $N_{n-k}(X)$ when we work with dimension instead of codimension. The direct sum $$N(X):=\bigoplus_{k=0}^nN^k(X)$$ is a graded $\mathbb R$-algebra with multiplication induced by the intersection form.
\par Define the cones $\Eff^i(X)=\Eff_{n-i}(X)$ as the closures of the cones of effective cycles in $N^i(X)$. The elements of $\Eff_i(X)$ are usually called \textit{pseudoeffective}. Dually, we have the nef cones $$\Neff^k(X):=\{\alpha\in\Eff^k(X)|\ \alpha\cdot\beta\geq 0\ \forall\beta\in\Eff_k(X)\}.$$  

\noindent Now let $C$ be a smooth complex projective curve and let $E$ be a locally free sheaf on $C$ of rank $n$, degree $d$ and slope $\mu(E):=\frac{\deg E}{\rk E}$, or $\mu$ for short. Let $$\pi:\mathbb P(E)\to C$$ be the associated projective bundle of quotients of $E$. The graded algebra $N(\bP(E))$ is generated in degree 1 by the classes $f$ and $\xi$ of a fiber of $\pi$ and of the Serre $\mathcal O_{\bP(E)}(1)$ sheaf respectively. It is completely described by:

\begin{equation}\label{i:e2} f^2=0,\qquad \xi^{n-1}f=[pt],\qquad \xi^n=d\cdot[pt].\end{equation}

\begin{figure}
\begin{center}
\includegraphics[scale=1]{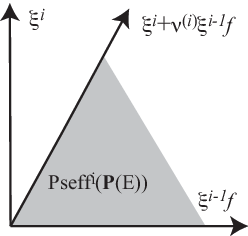}
\caption{$\Eff^i(\bP(E))$}
\label{Fig1}
\end{center}
\end{figure}

\noindent As a consequence of previous remarks, $N^i(\bP(E))$ and $\Eff^i(\bP(E))$ are 2-dimensional in positive dimension and codimension and $\xi^{i-1}f$ is a boundary of the later for $i\in\{1,\ldots,n-1\}$. The other boundary is spanned by $\xi^i+\nu^{(i)}\xi^{i-1}f$, which defines $\nu^{(i)}=\nu_{n-i}$. See Figure \ref{Fig1}. 

\noindent $E$ has a Harder--Narasimhan filtration (Prop 6.4.7, \cite{Lazarsfeld2}), $$0=E_l\subset E_{l-1}\subset\ldots\subset E_1\subset E_0=E,$$ for some $l$. Recall that by definition the successive quotients $Q_i:=E_{i-1}/E_i$ are semistable and their slopes $\mu_i:=\mu(Q_i)$ form an ascending sequence. The following theorem computes $\nu_i$ in terms of all the numerical data appearing in the Harder--Narasimhan filtration of $E$.

\begin{thm}\label{p:t} With the above setup, let $$r_i:=\rk Q_i,\ d_i:=\deg Q_i,\ \underline r_i:=\rk (E/E_i)=\sum_{k=1}^ir_k,\ \underline d_i:=\deg E_i=d-\sum_{j=1}^id_j.$$ Then, for all $k\in\{1,\ldots,l\}$ and $i\in\{1,\ldots,r_k\}$, except when $k=l$ and $i=r_l$: \begin{equation}\label{i:e'}\nu_{\underline r_{k-1}+i}=\nu^{(n-\underline r_{k-1}-i)}=-\underline d_{k-1}+i\mu_k\end{equation}
\end{thm}

\begin{figure}[ht]
\begin{center}
\includegraphics[scale=1]{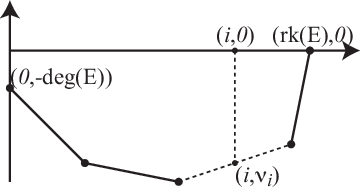}
\caption{Reading the boundaries of $\Eff_i(\bP(E))$ from $\mathcal P$}\label{Fig2}
\end{center}
\end{figure}

\noindent The formulas can be extracted from a picture strongly resembling the one in the Shatz stratification (Ch 11, \cite{LeP}). Construct the polygonal line $\mathcal P$ joining the points of coordinates $(\underline r_k,-\underline d_k)$ for $k\in\{1,\ldots,l\}$. See Figure \ref{Fig2}.
The theorem implies that the points of coordinates $(i,\nu_i)$ all lie on $\mathcal P$ for $i\in\{1,\ldots,n-1\}$. Note that the sides of this polygonal line have slopes $\mu_1,\ldots,\mu_l$ in this order.
\par Relating properties of objects on a projective bundle over a curve to the associated Harder--Narasimhan filtration is also apparent in work of H. Chen (\cite{Chen}) and A. Wolfe (\cite{Wolfe}) who independently computed the volume function on $\Eff^1(\bP(E))$ in terms of the numerical information of the Harder--Narasimhan filtration of $E$.  

\par For the proof, we start with the semistable case which is covered by a generalization of a result of Miyaoka (see \cite{Miyaoka}). We use the notation $\langle A\rangle$ for the convex cone spanned by a subset $A$ of a real vector space. 

\begin{prop}\label{p:ss} If $E$ is semistable of rank $n$ and slope $\mu$, then for all $i\in\{1,\ldots,n-1\}$, $$\Eff^i(\bP(E))=\left\langle (\xi-\mu  f)^i,\ \xi^{i-1}f\right\rangle.$$\end{prop} 

\noindent When $E$ is unstable, we have a natural inclusion as a proper subvariety $\imath:\bP(Q_1)\to \bP(E)$ and there is a fiber-by-fiber linear projection map $p:\bP(E)\setminus\bP(Q_1)\to\bP(E_1)$. We then perform induction showing that low dimensional cycles on $\bP(E)$ come from $\bP(Q_1)$, while higher dimensional cycles are related to cycles on $\bP(E_1)$ as illustrated by the following two assertions.

\begin{prop}\label{p:us1} The equality $\Eff_i(\bP(E))=\left\langle[\bP(Q_1)]\cdot(\xi-\mu_1f)^{r_1-i},\ \xi^{n-i-1}f\right\rangle$ holds for all $i\in\{1,\ldots,r_1\}$.\end{prop}
\noindent In fact, for $i<r_1$ the pushforward $\imath_*$ induces an isomorphism $\Eff_i(\bP(Q_1))\simeq\Eff_i(\bP(E))$ and for $i=r_1$, $$\Eff_{r_1}(\bP(E))=\left\langle[\bP(Q_1)],\ \xi^{n-r_1-1}f\right\rangle.$$
The projection $p:\bP(E)\setminus\bP(Q_1)\to\bP(E_1)$ induces for all $i$ a map $$\cone^{(i)}:\Eff^i(\bP(E_1))\to\Eff^i(\bP(E))$$ that geometrically sends a subvariety $Z\subset\bP(E_1)$ to the closure $\overline {p^{-1}(Z)}$ in $\bP(E)$ i.e. the cone over $Z$ with center $\bP(Q_1)$. The following proposition shows that every high dimensional cycle on $\bP(E)$ is equivalent to a cone over a cycle in $\bP(E_1)$.

\begin{prop}\label{p:us2} The map $\cone^{(i)}:\Eff^i(\bP(E_1))\to\Eff^i(\bP(E))$ is an isomorphism for $i\leq n-r_1-1$.\end{prop}
\noindent A rigorous construction for the coning map will be given in the proof where it will be clear why it is well defined. The proof of Proposition \ref{p:us2} is the more technical part of the main result.

\vskip.2cm \noindent The following statement is a consequence of Theorem \ref{p:t} and will be proved as Lemma \ref{efnef}:
\begin{prop} Let $C$ be a smooth projective curve and let $E$ be a locally free sheaf of rank $n$ on $C$. Then $E$ is semistable if, and only if, for all (any) $k\in\{1,\ldots,n-1\}$ we have that $\Eff^k(\bP(E))=\Neff^k(\bP(E))$.\end{prop}

It is natural to ask what happens if we work over an arbitrary smooth complex projective polarized variety $(V,H)$ with slope semistability in the sense of Mumford. Thomas Peternell suggests that if $E$ is slope unstable on $V$, then one should be able to find a pseudoeffective but not nef cycle on $\bP(E)$ in a natural way. In the application at the end of the paper we prove this result. We also construct an example showing that pseudo-effectivity and nefness need not be equivalent properties for cycles on $\bP(E)$ with $E$ a slope semistable bundle on $\mathbb P^2$.

\begin{paragraph}{Acknowledgments.} The author is greatly indebted to Robert Lazarsfeld for suggesting the main problem, and for sharing his intuition on many aspects of this paper. Thomas Peternell suggested the application, for which the author is grateful. Thanks also go to William Fulton and to Victor Lozovanu for helpful discussions.\end{paragraph}

\section{The proof of the main result}

We recall notation. $E$ is a locally free sheaf of degree $d$ and rank $n$ on a smooth complex projective curve $C$. Inside $N(\bP(E))$, $\xi$ and $f$ denote the classes of the Serre $\mathcal O(1)$ bundle on $\bP(E)$ and that of a fiber of the projection $\pi:\bP(E)\to C$ respectively. The sheaf $E$ admits a Harder--Narasimhan filtration  $E=E_0\supset E_1\supset\ldots\supset E_l=0$ and $Q_i:=E_{i-1}/E_i$. Denote $r_i=\rk(Q_i)$, $d_i=\deg(Q_i)$ and $\mu_i=\mu(Q_i):=\frac{d_i}{r_i}$. We set $X:=\bP(E)$ and start by describing the cone of nef divisors $\Nef(X)=\Neff^1(X)$.

\begin{lem}[Miyaoka]\label{pf:nef} $\Nef(X)=\langle \xi-\mu_1f,f\rangle$.\end{lem}
\begin{proof}
Hartshorne's Theorem (\cite{Hartshorne1}, or Thm 6.4.15, \cite{Lazarsfeld2}) states that the twist (in the sense of Section 6.3, \cite{Lazarsfeld2}) of $E$ by a $\mathbb Q$-divisor $\delta$ on $C$, $E\langle\delta\rangle$, is nef if and only if it has no negative slope quotient. Basic properties of the Harder--Narasimhan filtration say that $Q_1\langle\delta\rangle$ has the smallest slope among all quotients of $E\langle\delta\rangle$. Since $\mu(Q_1\langle\xi-\mu_1f\rangle)=0$, $\xi-\mu_1f$ generates a boundary of the nef cone. The other boundary is trivially spanned by $f$ and the result follows.\end{proof}

We are now ready to treat the semistable case, generalizing a result of Miyaoka (see \cite{Miyaoka}). 

\begin{lem}\label{pf:ss} If $E$ is semistable of rank $n$ and slope $\mu$, then for all $i\in\{1,\ldots,n-1\}$, $$\Eff^i(X)=\left\langle (\xi-\mu  f)^i,\ \xi^{i-1}f\right\rangle.$$\end{lem}
\begin{proof} Because there is only one term in the Harder--Narasimhan filtration of $E$, we have $\mu=\mu_1$ and by the intersection relations (\ref{i:e2}), $$(\xi-\mu_1f)^n=0.$$ From Lemma \ref{pf:nef}, $(\xi-\mu f)^i$ and $(\xi-\mu f)^{i-1}f=\xi^{i-1}f$ are intersections of nef divisors, so they are pseudoeffective. Conversely, if $a(\xi-\mu f)^i+b\xi^{i-1}f$ is pseudoeffective, then intersecting with $(\xi-\mu f)^{n-i}$ and $\xi^{n-i-1}f$ and using previous remarks shows that $a$ and $b$ are non-negative.\end{proof}

\par Our main effort is to study the case when $E$ is unstable. Assuming this, let $$0\to E_1\to E\to Q_1\to 0$$ be the short exact sequence induced by the Harder--Narasimhan filtration of $E$ with $Q_1$ the largest most negative slope quotient of $E$. Recall that $\imath:\bP(Q_1)\to \bP(E)$ denotes the canonical embedding. A slight generalization of Lemma \ref{pf:ss} allows to tie cycles of dimension at most $r_1$ on $\mathbb P(E)$ to cycles on $\mathbb P(Q_1)$.

\begin{lem}\label{pf:us1} For all $i\in\{1,\ldots,r_1\}$, $$\Eff_i(X)=\left\langle[\bP(Q_1)]\cdot(\xi-\mu_1f)^{r_1-i},\ \xi^{n-i-1}f\right\rangle.$$ In particular,  $\imath_*$ induces an isomorphism $\Eff_i(\bP(Q_1))\simeq\Eff_i(X)$ for $i<r_1$.\end{lem}
\begin{proof} The result in (Ex 3.2.17, \cite{Fulton}) adjusted to bundles of quotients over curves shows that $$[\mathbb P(Q_1)]=\xi^{n-r_1}+(d_1-d)\xi^{n-r_1-1}f.$$ Since $\xi-\mu_1f$ is nef, 
$$\tau_i:=[\mathbb P(Q_1)]\cdot(\xi-\mu_1f)^{r_1-i}=(\xi^{n-r_1}+(d_1-d)\xi^{n-r_1-1}f)(\xi-\mu_1f)^{r_1-i}\in \Eff_i(X)$$ for all $i\in\{1,\ldots,r_1\}$. Because $\{1,\xi,\ldots,\xi^{n-1}\}$ are linearly independent over $N(C)$, $\tau_i$ is nonzero for $i\in\{1,\ldots,r_1\}$. To see that they are actually in the boundary, use the nefness of $\xi-\mu_1f$ and notice that \begin{equation}\label{p:e1}\tau_i\cdot(\xi-\mu_1f)^i=(\xi^{n-r_1}+(d_1-d)\xi^{n-r_1-1}f)(\xi^{r_1}-r_1\mu_1\xi^{r_1-1}f)=0.\end{equation} 
That $\imath_*$ induces an isomorphism between the pseudoeffective cones follows from Lemma \ref{pf:ss} for the semistable bundle $Q_1$, from $\imath^*\mathcal O_{\bP(E)}(1)=\mathcal O_{\bP(Q_1)}(1)$ and from the projection formula. 
\end{proof}
Note that $$\tau_i=\xi^{n-i}+(d_1-d-\mu_1(r_1-i))\xi^{n-i-1}f=\xi^{n-i}+(-d+i\mu_1)\xi^{n-i-1}f,$$ so Theorem \ref{p:t} is proved for $k=1$ and $i\in\{1,\ldots,r_1\}$.

We move on to describe the coning construction that will allow us to tie the cycles not covered by the previous lemma to cycles on $\bP(E_1)$. Let $Y:=\mathbb P(E_1)$ and let $\rho:Y\to C$ be its bundle map. The projection map $$p:\bP(E)\setminus\bP(Q_1)\to \bP(E_1)$$ can be seen as a rational map $\xymatrix{X\ar@{.>}[r]& Y}$ whose indeterminacies are resolved by blowing--up $\bP(Q_1)$. Denote this blow-up by $\widetilde X$. There is a commutative diagram:
\begin{equation}\label{sq1}\xymatrix{\widetilde X=\mbox{Bl}_{\bP(Q_1)}\bP(E)\ar[r]^{\eta}\ar[d]^{B}& \mathbb P(E_1)=Y\ar[d]^{\rho}\\ X=\mathbb P(E)\ar[r]^{\pi}& C}\end{equation} where $B$ is the blow-down and $\eta$ is the resolved map mentioned in the above.   
\par The map $\cone^{(i)}$ is defined as the restriction of $B_*\eta^*$ to $\Eff^i(\bP(E_1))$. Before we can say anything about $\cone^{(i)}$, we need to know more about $\widetilde X$ and its intersection theory. This is achieved by the following proposition which describes $(\widetilde X,Y,\eta)$ as a projective bundle over $Y$ with fiber $\mathbb P^{r_1}$.

\begin{prop}\label{bup}$\mbox{}$ 
\begin{enumerate}
\item With the above notation, there exists naturally a locally free sheaf $F$ on $Y$ such that $\widetilde X\simeq \mathbb P_Y(F)$ and $\eta:\mathbb P_Y(F)\to Y$ is its associated bundle map. 

\item Let $\xi_1$ be the class of $\mathcal O_{\mathbb P(E_1)}(1)$, $f_1$ the class of a fiber of $\rho$, $\gamma$ the class of $\mathcal O_{\mathbb P_Y(F)}(1)$ and $\widetilde E$ the class of the exceptional divisor of $B$. We have the following change of bases relations:
\begin{equation}\label{b1}
\gamma=B^*\xi,\qquad
\eta^*\xi_1=B^*\xi-\widetilde E,\qquad 
\eta^*f_1=B^*f\end{equation}

\item The space $N(\widetilde X)$ is a free $N(Y)$-module via the pullback map $\eta^*$ and
\begin{equation}\label{b5}\widetilde E\cdot B^*(\xi-\mu_1f)^{r_1}=0.\end{equation}

\item If by abuse $\widetilde E$ also denotes the support of the exceptional divisor of $\widetilde X$, then with $\jmath:\widetilde E\to\widetilde X$ the canonical inclusion, $\widetilde E\cdot N(\widetilde X)=\jmath_*N(\widetilde E)$ as subsets of $N(\widetilde X)$. 
\end{enumerate}
\end{prop}

\begin{proof} The first line of the following commutative diagram induces the second, defining $F$:
$$\xymatrix{0\ar[r]& \rho^*E_1\ar@{->>}[d]\ar[r]&\rho^*E\ar[r]\ar[d]& \rho^*(Q_1)\ar[r]\ar@{=}[d]& 0
\\ 0\ar[r]& \mathcal O_{\mathbb P(E_1)}(1)\ar[r]& F\ar[r]& \rho^*(Q_1)\ar[r]& 0}
$$
\noindent Both lines are short exact sequences and the first vertical map is the tautological surjection. 
\par Let $X'=\mathbb P_Y(F)$ and $\lambda:X'\to Y$ be the bundle projection. From the Snake Lemma, $\rho^*E\to F$ is a surjective map and together with $\lambda^*F\to\mathcal O_{\mathbb P_Y(F)}(1)$ induces a surjective morphism $\lambda^*\rho^*E\to\mathcal O_{\mathbb P_Y(F)}(1)$ that determines $\sigma: X'\to X$ with $\sigma^*\mathcal O_X(1)=\mathcal O_{\mathbb P_Y(F)}(1)$. We want to show that we can identify $(\widetilde X,\eta,B)$ and $(X',\lambda,\sigma)$. We also have the commutative diagram:

\begin{equation}\label{sq2}\xymatrix{X'=\mathbb P_Y(F)\ar@/{}^{1pc}/[drr]^{\lambda}\ar@/{}_{1pc}/[ddr]_{\sigma}\ar[dr]^i\\
& X\times_C Y=\mathbb P_Y(\rho^* E)\ar[r]^-{\pr_2}\ar[d]^{\pr_1} & \mathbb P(E_1)=Y\ar[d]^{\rho}\\ 
& X=\mathbb P(E)\ar[r]^{\pi}& C}\end{equation}
In the above, $i$ is induced by the universality property of the fiber product and by the onto morphism $\rho^*E\to F$. In particular, $i$ is a closed immersion. 

\par The image of the composition $\pi^*E_1\to\pi^*E\to\mathcal O_{\mathbb P(E)}(1)$ is $\mathcal I\otimes\mathcal O_{\mathbb P(E)}(1)$, where $\mathcal I$ is the ideal sheaf of $\mathbb P(Q_1)$ in $X$. If $\mathcal S$ denotes the $\mathcal O_X$ algebra $\mathcal O_X\oplus\mathcal I\oplus \mathcal I^2\oplus\ldots$, then we have an induced surjective map of graded $\mathcal O_X$-algebras $\Sym(\pi^*E_1)\to \mathcal S*\mathcal O_X(1)$ with the notation in (\cite[Ch.II.7]{Hartshorne}: $\mathcal S*\mathcal L:=\oplus_{i\geq0}I^i\otimes \mathcal L^{\otimes i}$ for any invertible sheaf $\mathcal L$). This induces the closed immersion $i':\widetilde X=\Proj(\mathcal S*\mathcal O_{\mathbb P(E)}(1))\to \Proj(\Sym(\pi^*E_1))=X\times_CY$ that will fit inside a diagram similar to (\ref{sq2}). In particular $B$ and $\eta$ factor through $\pr_1$ and $\pr_2$. 
\par We have proved that $\widetilde X$ and $X'$ lie inside $X\times_CY$ and we want to prove that $(\widetilde X,B,\eta)$ and $(X',\sigma,\lambda)$ are equal. Since $\lambda$ and $\eta$ factor through $\pr_2$ while $\sigma$ and $B$ factor through $\pr_1$, it is enough to show that $X'=\widetilde X$. And because we are working over an algebraically closed field, it suffices to prove this over the closed points of $C$. Now the result is analogous to \cite[Ex.II.2.11.4]{Hartshorne}.

\par For the change of bases formulas in part $(ii)$, recall that $\sigma^*\mathcal O_X(1)=\mathcal O_{X'}(1)$ which yields $B^*\xi=\gamma$. That $B^*f=\eta^*f_1$ is a consequence of the commutativity of the square in diagram (\ref{sq1}). 
\par The closed immersion $i=i'$ in diagram (\ref{sq2}) induces a compatibility between the associated $\mathcal O(1)$ sheaves of $\Proj(\Sym(\pi^*E_1))=X\times_CY$ and $\Proj(\mathcal S*\mathcal O_X(1))=\widetilde X$. For the projective bundle $\Proj(\Sym(\pi^*E_1))$, this $\mathcal O(1)$ sheaf is $\pr_2^*\mathcal O_{\mathbb P(E_1)}(1)$. For $\Proj(\mathcal S*\mathcal O_X(1))$, the associated invertible sheaf is, by (Lem 7.9 and the proof of Prop 7.13, \cite{Hartshorne}), $\mathcal O_{\widetilde X}(-\widetilde E)\otimes B^*\mathcal O_{\mathbb P(E)}(1)$. Since $\eta$ factors through $\pr_2$, it follows that $\eta^*\mathcal O_{\mathbb P(E_1)}(1)=\mathcal O_{\widetilde X}(-\widetilde E)\otimes B^*\mathcal O_{\mathbb P(E)}(1)$ which yields $B^*\xi=\eta^*\xi_1+\widetilde E$.

\par The extension $0\to\mathcal O_Y(1)\to F\to \rho^*(Q_1)\to 0$ determines the total Chern class relation $$c(F)=c(\mathcal O_Y(1))\cdot c(\rho^*(Q_1))=(1+\xi_1)\cdot\rho^*(1+d_1\cdot [pt])=(1+\xi_1)(1+d_1f_1).$$
\noindent Plugging this into the appropriate Grothendieck relation and using (\ref{b1}), (\ref{b5}) follows easily.

\par The variety $\mathbb P(Q_1)\times_CY$ is the full preimage of $\mathbb P(Q_1)$ in $X\times_CY$ via $\pr_1$ and has the same dimension as $\widetilde E$ which shows that they are equal. To justify the equality $\widetilde E\cdot N(\widetilde X)=\jmath_*N(\widetilde E)$, one uses an explicit description of $N(\widetilde E)$ and $N(\widetilde X)$ as free modules over $N(Y)$ and the projection formula.
\end{proof}

\begin{rem}\label{bupr} The description of the blow-up as a projective bundle remains valid if $C$ is replaced by any nonsingular variety. The relations (\ref{b1}) remain true, but they no longer represent a change of bases.\end{rem}  

\begin{defn} If $V$ and $W$ are smooth varieties, we call a map $\varphi:N^i(V)\to N^i(W)$ pseudoeffective if $\varphi(\Eff^i(V))\subset \Eff^i(W)$.\end{defn}
We next relate pseudoeffective cycles of dimension bigger than $r_1=\rk Q_1$ on $X=\bP(E)$ to pseudoeffective cycles on $Y=\bP(E_1)$ using the coning construction. 

\begin{lem}\label{pf:us2} The map $\cone^{(i)}:=B_*\eta^*|_{\Eff^i(\bP(E_1))}$ is an isomorphism onto $\Eff^i(\bP(E))$ for $i<n-r_1$. \end{lem}
\begin{proof}
\par There is a visible isomorphism of abstract groups $\phi_i:N^i(X)\to N^i(Y)$ for $i<\dim Y=n-r_1$ sending $a\xi^i+b\xi^{i-1}f$ to $a\xi_1^i+b\xi^{i-1}_1f_1$. We prove that it induces an isomorphism $\Eff^i(X)\simeq \Eff^i(Y)$ for all such $i$, but for this we need more geometric descriptions for $\phi_i$ and its inverse. Define $U_i:N^i(Y)\to N^i(X)$ by 
$$U_i(c)=B_*\eta^*c$$ 
This is precisely the "coning" construction. $U_i$ is well defined since $\eta$ is flat and $B$ is birational. It is also clear that $U_i$ is pseudoeffective. We now check that $U_i=\phi_i^{-1}$. For this we will make use of the change of basis relations (\ref{b1}) and the projection formula.
$$U_i(a\xi_1^i+b\xi_1^{i-1}f_1)=B_*(a(B^*\xi-\widetilde E)^{i}+b(B^*\xi-\widetilde E)^{i-1}\cdot B^*f)$$ 
Expanding in the last formula shows that in excess of what we are looking for, there is a sum of the form $B_*(\sum_{1\leq j\leq i}\widetilde E^j\cdot B^*(\alpha_{i,j}))$ for some cycles $\alpha_{i,j}\in N(X)$ of varying dimensions. To show that this vanishes, it is enough, by the projection formula, to see that $B_*(\widetilde E^j)=0$ for all $j\leq i<n-r_1$. This is because $\widetilde E^j$ has dimension $n-j>r_1=\dim(\mathbb P(Q_1))$, so $B$ contracts it. Thus $U_i(a\xi_1^i+b\xi_1^{i-1}f_1)=a\xi^i+b\xi^{i-1}f$.

\par We construct an inverse for $U_i$ and prove that it is also pseudoeffective. Put $\delta=B^*(\xi-\mu_1f)^{r_1}$ and define $D_i:N^i(X)\to N^i(Y)$ by $$D_i(k)=\eta_*(\delta\cdot B^*k).$$ 
We show that $D_i=\phi_i$. By definition, $D_i(a\xi^i+b\xi^{i-1}f)=\eta_*(\delta\cdot B^*(a\xi^i+b\xi^{i-1}f)).$ Modulo $\widetilde E$, by (\ref{b1}), $B^*(a\xi^i+b\xi^{i-1}f)$ is $\eta^*(a\xi^i_1+b\xi^{i-1}_1f_1)$ and since $\delta\cdot\widetilde E=0$ by (\ref{b5}), one gets:
$$D_i(a\xi^i+b\xi^{i-1}f)=\eta_*(\eta^*(a\xi^i_1+b\xi^{i-1}_1f_1)\cdot\delta)=(a\xi^i_1+b\xi^{i-1}_1f_1)\cdot[Y]=a\xi^i_1+b\xi^{i-1}_1f_1.$$ 
We have used the projection formula and the identity $\eta_*\delta=[Y]$ which follows easily from (\ref{b1}) and $\eta_*\gamma^{r_1}=[Y]$. The later is a classical result (see Proof of Prop 3.1.a.i, \cite{Fulton}).
\par We still need to prove that $D_i$ is a  pseudoeffective map. For any effective cycle $k$ on $X$, $B^*k=k'+\jmath_*\widetilde k$, where $k'$ is an effective class (the strict transform under $B$), $\widetilde k$ is a not necessarily effective cycle class in $\widetilde E$ and $\jmath:\widetilde E\to\widetilde X$ is the canonical inclusion. Since $\delta$ is an intersection of nef classes and $\eta_*$ is pseudoeffective, it is enough to check that $\delta\cdot\jmath_*\widetilde k=0$ for any class in $\widetilde E$. This follows from (\ref{b5}) and the last part of Proposition \ref{bup}. The proof of the lemma is complete. \end{proof}

To finish the proof of Theorem \ref{p:t}, one applies induction noticing that the coning map is compatible with the most natural bases of $N(Y)$ and $N(X)$. We observe that deleting $E$ from its Harder--Narasimhan filtration amounts to deleting the first segment of the polygonal line $\mathcal P$ in Figure \ref{Fig2} if we assume by induction that the Theorem holds for $E_1$.
\par Also note that Lemma \ref{pf:us2} is vacuous when $E$ is semistable, or when $\rk (Q_1)=n-1$. However, Theorem \ref{p:t} is covered in these cases by Lemma \ref{pf:us1}.  


\section{An application}
The framework as well as the question to be answered by Proposition \ref{a1} were presented to the author by Thomas Peternell. Let $X$ be a smooth projective variety of dimension $n$ with a choice of an ample class $H$, and $E$ a locally free sheaf of rank $r$ on $X$. Recall that $\Eff^k(X)$ is the closed cone in $H^{2k}(X,\mathbb R)$ spanned by classes of codimension $k$ subvarieties of $X$, whereas $\Neff^k(X)$ is the closed cone spanned by classes in $\Eff^k(X)$ that have nonnegative intersection with members of $\Eff_k(X)$.

\begin{defn} We say $E$ is $k-$homogeneous if every pseudoeffective $k-$dimensional cycle on $\bP_X(E)$ is nef i.e. $\Eff^k(\bP(E))=\Neff^k(\bP(E))$.\end{defn} 
Recall that $E$ is slope semistable if for all nonzero coherent $F\subset E$, one has $\mu(F)\leq\mu(E)$ with $\mu(F):=\frac{c_1(F)\cdot H^{n-1}}{rk(F)}$ the $H-$slope of $F$. 

\begin{lem}\label{efnef} If $X$ is a curve, then a locally free sheaf $E$ of rank $r$ is semistable if, and only if, it is $k-$homogeneous for all (or for any) $k\in\{1,\ldots,r-1\}$.\end{lem}
\begin{proof} Up to a rational twist, one can assume that $\deg E=0$. If $E$ is semistable of degree $0$, then by Lemma \ref{pf:ss}, $\Eff^k(\bP(E))$ is spanned by $\xi^k$ and $\xi^{k-1}f$ for all $k\in\{1,\ldots,r-1\}$. Since $\xi^r=0$, it follows that $\chi:=a\xi^k+b\xi^{k-1}f$ is nef if, and only if, $$\chi\mbox{ is pseudoeffective, }b=\xi^{r-k}(a\xi^k+b\xi^{k-1}f)\geq 0, \mbox{ and }a=\xi^{r-k-1}f(a\xi^k+b\xi^{k-1}f)\geq 0.$$ This shows $\Neff^k(\bP(E))=\Eff^k(\bP(E))$ for all $k$.
\par Assume now that $E$ is unstable. Then Theorem \ref{p:t}, as illustrated in Figure \ref{Fig2}, proves that $\nu^{(k)}$ is negative for all $k\in\{1,\ldots,r-1\}$. Because $\deg E=0$, $\Nef^k(\bP(E))=\left\langle \xi^k-\nu_{k}\xi^{k-1}f,\ \xi^{k-1}f\right\rangle$. In particular, $\xi^k$ is effective but not nef and $E$ is not $k-$homogeneous for all $k\in\{1,\ldots,r-1\}$.\end{proof} 

\par The question is what happens if $X$ is of arbitrary dimension? More precisely, we prove:

\begin{prop}\label{a1} Assume that $E$ is slope--unstable with respect to $H$. Then there exists $k$ such that $E$ is not $k-$homogeneous.\end{prop}   
\begin{proof} We have to find some $k$ and two pseudoeffective cycles of codimension $k$ and dimension $k$ respectively on $\mathbb P(E)$ whose intersection is negative. The idea is to use the Mehta--Ramanathan theorem (see \cite{Mehta}) to restrict to the curve case, where we use Theorem \ref{p:t} to produce a $k-$dimensional pseudoeffective and not nef cycle whose pushforward to $\bP(E)$ we show enjoys the same properties.
\par The locally free sheaf $E$ admits a Harder--Narasimhan filtration $\ldots\subset E_1\subset E_0=E$ by torsion free subsheaves and let $Q=E/E_1$. The sheaf $Q$ is locally free of rank $s$ off a codimension at least 2 locus. We have a closed immersion $\mathbb P(Q)\to\mathbb P(E)$ and there is a unique irreducible component $Z$ of $\mathbb P(Q)$ that dominates $X$. We set $k=r-s$ and choose $[Z]\in N^k(\mathbb P(E))$ as our effective codimension $k$ cycle.
\par Let $C$ be a general complete intersection curve numerically equivalent to $(NH)^{n-1}$ for $N\gg0$. By rescaling $H$, one may assume $[C]=H^{n-1}$ in $N_1(X)$. It will be useful to assume that $\deg E|_C=0$ which we can by making a rational twist of $E$ by a multiple of $H$. By the Mehta--Ramanathan theorem, the Harder--Narasimhan filtration of $E$ restricts to the Harder--Narasimhan filtration of $E|_C$. From the assumptions that $E$ is not slope semistable and that $\deg E|_C=0$, it follows that $\deg Q|_C<0$. Consider the commutative diagram:

$$\xymatrix{(\bP_C(Q|_C),\zeta_C)\ar[ddr]_{\pi'}\ar[dr]^{i_C}\ar[r]^q & (\bP_X(Q),\zeta)\ar[dr]^i\ar@{.>}[ddr]& \\
& (\bP_C(E|_C),\xi_C)\ar[r]_(.45){e}\ar[d]^{\pi_C} & (\bP_X(E),\xi)\ar[d]^{\pi}\\ 
& C\ar[r]^j& X}$$
The maps $j,i,e,q,i_C$ are natural closed immersions. For example, $i$ is induced by $E\to Q$ and $e$ is induced by $E\to E|_C$. The numerical Serre $\mathcal O(1)$ classes $\xi,\zeta,\xi_C,\zeta_C$ on their respective projective bundles are compatible with the maps in the diagram i.e. $\zeta=i^*\xi$ etc. The morphisms $\pi, \pi'$ and $\pi_C$ are bundle projections. 
\par The cycle we are looking for is $\alpha=\xi^s\cdot \pi^*H^{n-1}$. The dimension of $\alpha$ is $n+r-1-(s+n-1)=r-s=k$. We need to show that $\alpha$ is pseudoeffective and that $\alpha\cdot[Z]<0$. Recall that $[C]=H^{n-1}$. 
$$\alpha\cdot[Z]=\xi^s\cdot\pi^*[C]\cdot i_*[\bP(Q)]$$ This is because we can choose $C$ so that no component of $\bP(Q)$ other than $Z$ meets its preimage. By the projection formula, $$\xi^s\cdot\pi^*[C]\cdot i_*[\bP(Q)]=i_*i^*(\xi^s\cdot\pi^*[C])=i_*(\zeta^s\cdot (\pi i)^*[C])=i_*(\zeta^s\cdot (\pi i)^*j_*[C]).$$ By base change and then again by the projection formula and the commutativity of the diagram above, the later is $$i_*(\zeta^s\cdot q_*[\bP(Q|_C)])=i_*q_*q^*\zeta^s=i_*q_*\zeta_C^s=e_*i_{C*}\zeta_C^s.$$ From the Grothendieck relation for $Q|_C$, we obtain that $\zeta_C^s=\deg Q|_C\cdot[pt]$. Since the degree of $Q|_C$ is negative and $e$ and $i_C$ are closed immersions, it follows that indeed $[Z]\cdot\alpha<0$.
\par We still need to show that $\alpha$ is pseudoeffective. Toward the end of the proof of Lemma \ref{efnef}, we have shown that $\xi_C^s$ is (pseudo)effective, hence the pushforward  $$e_*(\xi_C^s)=e_*e^*(\xi^s)=\xi^s\cdot e_*([\bP(E_C)])=\xi^s\cdot\pi^*[C]=\alpha$$ is also pseudoeffective. 
\end{proof}

\par The converse of Proposition \ref{a1} is in general false.
\begin{exap} There exists a rank 2 vector bundle on $\mathbb P^2$ sitting in an extension $$0\to\mathcal O_{\mathbb P^2}\to E\to J(1)\to 0,$$ where $J$ is the ideal sheaf of two distinct points in $\mathbb P^2$. Any such $E$ is stable, but not 1-homogeneous.\end{exap}
\begin{proof}The construction and stability of $E$ are explained in (Example 1, pag 187, \cite{Okonek}). If $\xi$ is the numerical class of $\mathcal O_{\mathbb P(E)}(1)$ on $\mathbb P(E)$, we show that $\xi$ is effective but not nef. The first assertion holds because $E$ has an obvious nonzero section. 
\par Let $\sigma:X\to\mathbb P^2$ be the blow-up of $\mathbb P^2$ along $J$, let $F$ be the exceptional divisor and $H$ the class of a line in $\mathbb P^2$. We have an epimorphism $\sigma^*E\to\mathcal O_X(-F)\otimes\sigma^*\mathcal O(1)$. The self intersection $(\sigma^*H-F)^2=\sigma^*H^2-2\sigma^*H\cdot F+F^2=1-0-2=-1$ is negative by the projection formula and because $F$ is the union of two disjoint $-1$ curves. Therefore $\mathcal O_X(-F)\otimes\sigma^*\mathcal O(1)$ is not nef showing that $\sigma^*E$, so $E$ and finally $\xi$ cannot be nef either.  
\end{proof}

\noindent\textsc{Department of Mathematics, University of Michigan, Ann Arbor, MI 48109, USA}
\textsc{E-mail:} mfulger@umich.edu
\vskip.5cm

\noindent\textsc{Institute of Mathematics of the Romanian Academy, P. O. Box 1-764, RO-014700,
Bucharest, Romania}

\end{document}